\documentclass{amsart}
\usepackage[leqno]{amsmath} 
\usepackage{amsfonts}%
\usepackage{amssymb}%
\usepackage{amsthm}%
\usepackage{verbatim}

\theoremstyle{plain}
\newtheorem{theo}{Theorem}[section]
\newtheorem{thm}{Theorem}[section]

\newtheorem{lema}[theo]{Lemma}

\newtheorem{coro}[theo]{Corollary}
\theoremstyle{definition}
\newtheorem{defi}[theo]{Definition}

\theoremstyle{remark}

\newcommand{\atiez}{\mbox{\quad and \quad}}

\begin{document}

\title{Constructive decomposition of functions of two variables using functions of one variable}

\author{Eva Trenklerov\'a}

\address{Department of Computer Science, Faculty of Science, 
P.J. \v{S}af\'arik University, Jesenn\'a 5, 04001 Ko\v{s}ice, Slovakia}
\email{eva.trenklerova@upjs.sk}
\thanks{The author was supported by grants VEGA 1/3002/06 and VEGA 1/3128/06}



\date{February 26, 2007}



\begin{abstract}

 Given a compact set $K$ in the plane, which contains no triple of points forming a vertical
and a horizontal segment, and a map $f\in C(K)$,
we give a construction of functions $g,h\in C(\mathbb R)$ such that
$f(x,y)=g(x)+h(y)$ for all $(x,y)\in K$. This provides a constructive proof of a part
of Sternfeld's theorem on basic embeddings in the plane.
In the  proof we construct a sequence of finite graphs,
which provide an arbitrarily good approximation of the set $K$.

\end{abstract}

\maketitle

\pagestyle{plain}



\newcommand{\cesta}[4]{{#1}_{#2} {#1}_{#3}\ldots {#1}_{#4}}
\newcommand{\Graph}[1]{\Gamma^{#1}}
\newcommand{\g}[1]{g^{#1}}
\newcommand{\f}[1]{f^{#1}}
\def\:{\mkern1mu\colon}
\newcommand{\fun}[3]{\mbox{${#1} \: {#2} \to {#3}$}}
\newcommand{\GraphPlus}[1]{\Gamma^{#1}_+}
\newcommand{\Gplus}[1]{\Gamma^{#1}_+}
\newcommand{\GplusNew}[1]{H^{#1}_{+}}
\newcommand{\GminusNew}[1]{H^{#1}_{-}}
\newcommand{\GraphMinus}[1]{\Gamma^{#1}_-}
\newcommand{\Gminus}[1]{\Gamma^{#1}_-}

\newcommand{\Gvert}[1]{\Gamma^{#1}_{\text{vert}}}
\newcommand{\Ghor}[1]{\Gamma^{#1}_{\text{hor}}}
\newcommand{\Gshort}[1]{\Gamma^{#1}_{\text{short}}}

\newcommand{\Vvert}{V_{\text{vert}}}
\newcommand{\Vhor}{V_{\text{hor}}}
\newcommand{\Vshort}{V_{\text{short}}}
\newcommand{\Vplus}{V_{+}}
\newcommand{\Vminus}{V_{-}}
\newcommand{\Rset}{\mathbb R}
\newcommand{\Nset}{\mathbb N}
\newcommand{\card}{\mathrm{Card}}

\section{Introduction}

An embedding (continous one-to-one map with  continuous inverse) $\fun \varphi K {\mathbb R^k}$
of a compactum (compact metric space) $K$ into the $k$-dimensional Euclidean space $\mathbb R^k$
is called a {\em basic embedding\/} provided that for each continuous real-valued function $f\in C(K)$, there exist
 continuous real-valued functions of one real variable  $g_1,\ldots,g_k \in C(\mathbb R)$ such that
$f(x_1,\ldots,x_{k}) = g_1(x_1)+\ldots+g_{k}(x_{k})$, for all points
$(x_1,\ldots,x_{k})\in \varphi(K)$. We also say, that the set $\varphi(K)$ is {\em basically embedded\/} in 
$\mathbb R^k$.

The question of the existence of basic embeddings was already implicitly
contained in Hilbert's 13th problem: Hilbert conjectured, that not all continuous functions
of three variables were expressible as sums and superpositions of continuous functions
of a smaller number of variables.

Ostrand~\cite{ostrand} proved, that each $n$-dimensional compactum 
can be basically embedded into $\mathbb R^{2n+1}$, for $n\geq 0$.
His result is a generalization of results by Arnold~\cite{ArnoldHilbert, ArnoldHilbertOne} and Kolmogorov~\cite{kolmogorovOne,kolmogorov}.
Sternfeld~\cite{sternfeld1} proved, that the parameter $2n+1$ is the best possible in a very strong sense:
namely, that no $n$-dimensional compactum can be basically embedded into $\mathbb R^{2n}$, for $n\geq 2$.
Ostrand's and Sternfeld's results thus characterize compacta basically embeddeble into $\mathbb R^k$ for $k\geq 3$.
Basic embeddability into the real line is trivially equivalent to embeddability.
The remaining problem of the chraracterization of compacta basically embeddable into
  $\mathbb R^2$, which was already raised by Arnold~\cite{ArnoldProblemSix},
was solved by Sternfeld~\cite{sternfeld2}:

\begin{thm}[Sternfeld]\label{v:SternfeldPlane}
Let $K$ be a compactum and let $\fun \varphi K {\mathbb R^k}$ be an embedding. Then

(B) $\varphi$ is a basic embedding if and only if

(A) there exists an $m\in \mathbb N$ such that the set $\varphi(K)$ contians no
array on $m$ points.
\end{thm}

\begin{defi} An {\em array \/} is a sequence of points in the plane $\{a_i\}_{i\in I}$,
where $I=\{1,2,\ldots,m\}$ or $I=\mathbb N$, such that for each $i$:
\begin{itemize}
\item $a_i\ne a_{i+1}$ and $[a_i;a_{i+1}]$ is a segment parallel to one of the coordinate axes and
\item the segments $[a_i;a_{i+1}]$ and $[a_{i+1};a_{i+2}]$ are mutually orthogonal.
\end{itemize}
\end{defi}

Using the geometric description (A),  Skopenkov~\cite{skopenkov} gave a characterization of continua basically embeddable into the 
plane by means of forbidden subsets; the characterization
resembles Kuratowski's characterization of planar graphs. 
In a similar way Kurlin~\cite{Kurlin00} characterized finite graphs
basically embeddable into $\mathbb R\times T_n$, where $T_n$ is
a star with $n$-rays.
Repov\v s and \v Zeljko~\cite{Zeljko} proved a result concerning the smoothness of functions
in a basic embedding in the plane.

However, Sternfeld's proof of  the equivalence (A) $\Leftrightarrow$  (B) 
is not direct,  but
uses a reduction to linear operators. In particular it is not constructive. 
It is therefore desirable to find a straightforward, constructive proof
which will consequently provide an elementary proof
of Skopenkov's and Kurlin's characterizations.
%
A constructive proof of  (B) $\Rightarrow$  (A)  is given
in~\cite{EvaNeza}. 

In this paper we give an elementary constructive proof of the implication (A) $\Rightarrow$ (B),
provided that $m=3$:

\begin{thm}\label{v:Basic} Let $\fun \varphi K {\mathbb R^2}$ be an embedding of a compactum
$K$ into the plane such that the set $\varphi(K)$ contains no array on three points.  
Then for every function $f \in C(\varphi(K))$
 there exist functions  $g,h\in C(\mathbb R)$ such that $f(x,y)=g(x)+h(y)$
for all points $(x,y)\in \varphi(K)$.
\end{thm}

It turns out,  that even if  the set $\varphi(K)\subseteq \mathbb R^2$ contains no array on
three points, 
a constructive decomposition of a function $f\in C(\varphi(K))$ is a non-trivial problem.
We also believe, that the proof can be modified as so as to obtain a constructive proof of
the implication   (A) $\Rightarrow$ (B) for an arbitrary $m \in \mathbb N$.

The author would like to thank an anonymous referee for pointing out a short and elegant proof of
Theorem~\ref{v:Basic} which  however  does not give an explicit construction.
We give his proof in the last section.

%

The idea of our proof is to approximate the set $\varphi(K)$ by a finite sequence graphs $\{\Graph n\}_{n\in \mathbb N}$;
with incresing $n$ the approximation becomes finer.
The property that $\varphi(K)$ contians no array on three points
 ensures that if the approximation is fine enough, then the graph
 does not contain any subgraph which \lq\lq resembles\rq\rq\ an array on three points
(Theorem~\ref{v:Veta1}).
For a given $f\in C(\varphi(K))$ and $\varepsilon>0$ we consider a fine-enough graph $\Graph n$
and  construct a function $\fun {\g n} {V(\Graph n)} {\mathbb R}$  on the vertices of the graph (Theorem~\ref{v:Veta2}).
Using $\g n$ we  construct functions $g,h \in C(\mathbb R)$
whose sum on $\varphi(K)$ does not differ from $f$ by more than a constant times $\varepsilon$ 
and whose norms are bounded by $2||f||$ (Theorem~\ref{v:Veta3}). By an iterative procedure we  obtain 
a decomposition of $f$ on $\varphi(K)$ (Theorem~\ref{v:Veta4}).

The author would like to than Du\v san Repov\v s and Arkadyi Skopenkov for the inspiration for
this paper, Lev Bukovsk\'y for support and especially Ne\v za Mramor-Kosta for endless conversations on the topic and
 invaluable advice.
\section{Proof}

Throughout the text we fix an embedding $\fun \varphi K {\mathbb R^2}$ of a compactum
$K$ into the plane such that the set $\varphi(K)$ contains no array on three points. 
For simplicity of notation we identify the set $K$ and its homeomorphic image $\varphi(K)$
and we shall speak about a set $K\subseteq \mathbb R^2$.
We also
fix a function $f\in C(K)$, a positive real $\varepsilon>0$
and a positive real $\delta>0$ such that 
\begin{equation}\label{e:Veta2ChoiceDelta}|x_1-x_2|<2\delta \text{\ \ implies that\
\ } |f(x_1)-f(x_2)|<\varepsilon \end{equation} for all $x_1,x_2\in K.$
We denote by $\fun {p,q} {\mathbb R^2} {\mathbb R}$ the vertical and horizontal orthogonal projections:
$p(x,y)=x$, $q(x,y)=y$. 
We  denote points in $\mathbb R^2$ by $(x,y)$,  we 
denote intervals in $\mathbb R$ by $(x;y), [x;y]$, etc. and segments by $[x;y]$.

For every $n\in \mathbb N$ we define  a unique graph $\Graph
n = (V(\Graph n ),E(\Graph n))$,  which reflects certain
properties of the set $K$:
Consider a point-lattice given by points
$(i/2^n,j/2^n)$
%
%
with $i,j\in \mathbb Z.$
 The vertex set $V(\Graph n)$ consists of  all points
$\left( {i}/{2^n},{j}/{2^n} \right)$  such that
the square
$\left[ {i}/{2^n};{(i+1)}/{2^n} \right) \times \left[ {j}/{2^n};{(j+1)}/{2^n} \right)$
meets the  set $K$.
The edge set $E(\Graph n)$ consists of all two-element sets $\{u_1,
u_2\}$ of vertices which have the same vertical or horizontal projection, or are neighbors
in the lattice in one of the projections, i.e.: 
$|p(u_1)-p(u_2)|\leq {1}/{2^n}$ or 
$|q(u_1)-q(u_2)|\leq {1}/{2^n}.$

%
%
An edge $\{u_1, u_2\}$ is denoted by $u_1u_2.$ If
$|p(u_1)-p(u_2)|\leq {1}/{2^n}$ then we say that the edge $u_1 u_2$
is {\em vertical\/}. If $|q(u_1)-q(u_2)|\leq {1}/{2^n}$ then we say
that the edge $u_1 u_2$ is {\em horizontal\/}. Hence, an edge may be
both vertical and horizontal.
A {\em path\/} is a sequence  $u_0, u_1,\ldots,u_k$ of pairwise
different vertices such that each two consecutive ones form an edge;
it is denoted by  $\cesta u 0 1   k$.

The following is the crucial statement about the geometric properties of the constructed graphs.
The result follows from the property that the set $K$ contains no array on three points.

\begin{theo}\label{v:Veta1}
For each $l\in \mathbb N$ and $\alpha>0$
there exists $n_0\in \mathbb N$ such that for all $n\geq n_0$ the following holds:
if $\cesta u 0 1   {k}$ is a path in $\Graph  n$ such that
\begin{equation}\label{e:Veta1}|p(u_0)-p(u_1)|\geq \alpha \atiez |q(u_{k-1})-q(u_{k})|\geq \alpha \end{equation} 
then its length is greater than $l$, i.e.: $$k > l.$$
\end{theo}

\begin{proof}
First, let us note the following. Let
$\{\Graph n \}_{n\in \mathbb N}$ be a sequence of graphs and
$\{u_1^{n} u_2^{n} \}_{n\in \mathbb N}$ a sequence of edges, with each
$u_1^{n} u_2^{n}\in E(\Graph n)$. 
Then, since $K$ is compact, there is a subsequence $\{m_n\}_{n\in \mathbb N}$
such that both $u_1^{m_n} \to u_1 \in K$
and $u_2^{m_n} \to u_2 \in K$ as $n\to \infty.$
From the graph construction  it follows, that either $u_1=u_2$ or $[u_1;u_2]$ is a segment
parallel to one of the coordinate axes.
Moreover, if $|p(u^{m_n}_1)-p(u^{m_n}_2)|\geq \alpha>0 $
for each $n$, then  $[u_1;u_2]$ is a segment 
parallel to the $x$-axis, with $|p(u_1)-p(u_2)|\geq \alpha$.
Similarly for the $q$-projection.

The proof of the theorem is indirect.
Assume, that for some $l_0$ and $\alpha_0$ there exists an infinite sequence of graphs
$\{\Graph {m_n}\}_{n\in \mathbb N}$ 
such that each graph
$\Graph {m_n}$ contains a path $u^{m_n}_0 u^{m_n}_1\ldots u^{m_n}_{k^{m_n}}$
for which (\ref{e:Veta1}) holds true but its length $k^{m_n}$ is not greater than $l_0$,
i.e. $k^{m_n}\leq l_0$.
Since the lentghs of the paths are limited by $l_0$, an infinite number of them has the same length.
Without loss of generality we may assume that $k^{m_n}= l_0$ for all $n$.

As we  noted in the beginning of the proof it follows, that there exist limit points 
$u_0,u_1,\ldots,u_{l_0} \in K$, such that  either $u_i=u_{i+1}$ or
$[u_i;u_{i+1}]$ is a segment parallel to one of the coordinate axes, for each $i$. In particular,
 (\ref{e:Veta1}) implies, that
$[u_0;u_1]$ is a horizontal segment and $[u_{l_0-1};u_{l_0}]$ is a vertical segment.
It follows, that the set $\{u_0,u_1,\ldots,u_{l_0}\}\subseteq K$ contains an array on at least three points.
%
%
%
\end{proof}

We call an edge $u_1 u_2\in E(\Graph n)$
{\em long\/} if $|u_1-u_2|\geq \delta$ and we call it
{\em short\/} if $|u_1-u_2|< \delta$.
Let
\begin{itemize}
\item $\Gshort n$ be the subgraph of $\Graph n$ induced by all short edges 
\item $\Ghor n$ be the subgraph of $\Graph n$ induced by all long horizontal edges 
\item $\Gvert n$ be the subgraph of $\Graph n$ induced by all long vertical edges. 
\end{itemize}

Theorem~\ref{v:Veta1} in particular implies that for a fine enough approximation the graph $\Graph n$
does not conatin any triple of vertices $u_1,u_2,u_3$ with $u_1u_2$ a long horizontal edge and $u_2u_3$
a long vertical edge. Such a triple \lq\lq resembles\rq\rq\ in a way an array on three points. 
Stated differently:
 
 \begin{coro} There exists $n_0\in \mathbb N$ such that for all $n\geq n_0$
the graphs  $\Ghor n$ and  $\Gvert n$ are disjoint.
\end{coro}

For each $n\in \mathbb N$ we fix a function
 $\fun {f^n}
{V(\Graph n)} {\mathbb R}$,
whose value in a given  vertex is equal to the value
of $f$ in some point of $K$ which lies in the corresponding square.
That is, for  each vertex $u\in V(\Graph n)$ we fix
a point $x_u$ which lies in
$K\cap ( \left[p(u);p(u) + {1}/{2^n}\right) \times
\left[q(u);q(u) + {1}/{2^n}\right)) $ and 
let $$f^n(u)=f(x_u).$$
%
%
%
Assumption (\ref{e:Veta2ChoiceDelta}) implies the following:

\begin{lema}\label{v:ZmenaFNaKratkej} Let $1/2^n<\delta$. If  $u_1u_2 \in E(\Gshort n)$ then $|\f n(u_1)-\f
 n(u_2)|<\varepsilon.$
\end{lema}

\begin{theo}\label{v:Veta2} 
There exists $n_0\in \mathbb N$ such that for all $n\geq n_0$ there exists a function
 $\fun {g^n}
{V(\Graph n)} {\mathbb R}$
  	 with the following properties:

\begin{enumerate}
\item For each edge $u_1u_2\in E(\Graph n)$:
\begin{enumerate}
\item if $u_1u_2\in E(\Gshort n)$  then $|g^n(u_1)-g^n (u_2)|\leq \varepsilon$
\item if $u_1u_2\in E(\Ghor n)$  then $|f^n(u_1)-g^n(u_1)|\leq \varepsilon$
\item if $u_1u_2\in E(\Gvert n)$  then $g^n(u_1)= 0$
\end{enumerate}
  \item $||g^n||\leq ||f^n||.$   
\end{enumerate}
\end{theo}

\begin{proof}
Denote
\begin{equation}\label{e:ChoiceOfL}
F=\left[\frac{||f||}{\varepsilon}\right]. \end{equation} 
Take $n_0$ from Theorem~\ref{v:Veta1} corresponding to
$l=F+1$ and $\alpha=\delta$. Take an $n\in \mathbb N$ such that
$n\geq n_0$ and $1/2^n<\delta$.


Let  $\Gplus n$ be the subgraph of $\Graph n$ induced  by the vertices $u$ with
$ \f n(u)\geq 0$ and let
$\Gminus n$ be the subgraph of $\Graph n$ induced by the vertices $u$ with
$ \f n(u)< 0$.
We shall define $\g n$  separately on $\Gplus n$ and $\Gminus n$.
So, consider $\Gplus n$.

If $\Gplus n \cap \Ghor n=\emptyset$ (i.e. the graph $\Gplus n$ contains no vertex,
which is the end of a long horizontal edge), then  let
$\g n(u)=0$ for all $u\in V(\Gplus n).$

If $\Gplus n \cap \Ghor n\not =\emptyset$, then we first modify the graph $\Gplus n$
by adding certain vertices and edges and removing some edges, and afterward we define $\g n$
on the modified graph. The modified graph is denoted by $\GplusNew n$.

Let
$$V(\GplusNew n) = V(\Gplus n) \cup \{w_0\ldots w_F\}.$$
%
We define $\f n$ on $ \{w_0\ldots w_F\}$ by letting
$$\f n (w_i) = i \varepsilon$$ for each $i$. Note that
$\f n(w_F) \approx ||f||$.
The edge-set $E(\GplusNew n)$ consists of
\begin{itemize}
\item the short edges from $E(\Gplus n)$ 
\item the edges $w_iw_{i+1}$ for all $i$
\item edges connecting some vertices of $V(\Gplus n)$ to the vertices $w_i$:
each vertex $u \in V(\Gplus n)\cap V(\Ghor n)$
is connected
to the vertex $w_i$ such that $[\f n(u)/ \varepsilon] = i$.
\end{itemize}
So, let
\begin{eqnarray*}
E(\GplusNew n)&=& E(\Gplus n \cap \Gshort n) \cup 
 \{w_0w_1,w_1w_2,\ldots, w_{F-1} w_{F} \} \cup \\
&&\bigcup_{u \in V(\Gplus n)\cap V(\Ghor n)}
  \{w_i u\mid  [\f n(u)/ \varepsilon] = i \}.
\end{eqnarray*}
The above definitions together with Lemma~\ref{v:ZmenaFNaKratkej} imply that
\begin{equation}\label{e:EachEdge} |\f n(u_1) - \f n (u_2)|\leq
\varepsilon
 \end{equation} for each edge $u_1u_2 \in E(\GplusNew n)$.

Now we define $\g n$ on $V(\GplusNew n)$.
%
%
Let $\fun d {V(\GplusNew n)} \{0,1,2,\ldots\}$ be the
function, which assigns to each vertex connected by a path to the vertex $w_F$
its distance from $w_F$, and which assigns  to each vertex not connected by a path
to $w_F$ the value 0.
 For each vertex $u\in
V(\GplusNew n)$ let
\begin{equation}\label{e:ThirdPartDefG} \g n(u) = \max\left\{
\left(F   - d(u)\right) \varepsilon,0
      \right\}. \end{equation}
Analogously we construct the graph $\GminusNew n$ and the define function $\fun {\g n} {V(\GminusNew n)} {\mathbb R}$.

Let us show that the function $\g n$ defined in this way satisfies all the
points from the statement of the theorem.

{\bf 1(a)} Let $u_1u_2\in E(\Gshort n)$.
If $u_1u_2\in E(\Gplus n)$ or $u_1u_2\in E(\Gminus n)$, then
 $|\g n(u_1)-\g n (u_2)|\leq \varepsilon$ follows directly from definition of $\g n$.
So, let 
$u_1\in V(\Gplus n)$ and $u_2\in V(\Gminus n)$. Using
(\ref{e:EachEdge}), (\ref{e:ThirdPartDefG}) and their analogies for $\GminusNew n$, by induction
on the distance from the vertex $w_F$, we can show that

\begin{equation}\label{e:WhereIsG} 
\begin{array}{l}
\g n(u) \in [0,\f n(u)],\  \forall u\in V(\GplusNew n) \\
\g n(u) \in [\f n(u),0], \  \forall u\in V(\GminusNew n) .
\end{array}
\end{equation}
So, by Lemma~\ref{v:ZmenaFNaKratkej} we have
$|\g n(u_1)-\g n (u_2)| = \g n(u_1)-\g n (u_2) \leq \f n(u_1)-\f n (u_2)\leq \varepsilon$.

%

{\bf 1(b)} Let $u_1u_2\in E(\Ghor n)$. Let for instance  $u_1\in V(\Gplus n)$.
We want to show that $|f^n(u_1)-g^n(u_1)|\leq \varepsilon$.
Since $u_1 \in V(\Ghor n)\cap V(\Gplus n)$, the edge-set $E(\GplusNew n)$
contains the edge $u_1w_i$ such that
$[\f n(u_1)/\varepsilon]=i$. 
%
 Therefore $\GplusNew n$ contains
the path $w_{F} w_{F-1} \ldots w_i
u$, so  $d(u_1)\leq F -i + 1.$
%
Hence $\g n(u_1) =(F-d(u_1)) \varepsilon\geq (i-1)\varepsilon  \geq \f n(u_1) -\varepsilon.$ 
On the other hand, equation (\ref{e:WhereIsG}) implies that $\g n(u_1)\leq \f n(u_1)$
and the statement follows.
 
%
%

{\bf 1(c)}
Let $u_1u_2\in E(\Gvert n)$. Let for instance  $u_1\in V(\Gplus n)$.
We want to show that $\g n(u_1)=0$. 
If $u_1$ is not connected by a path to $w_F$, then by definition $\g n(u_1)=0.$
 
So, let now
$P=w_F \ldots u_1$ be a path in $\GplusNew n$ such that $d(u_1)$ is
equal to its length. This path necessarily contains an edge of the form $w_iv$,
where $v\in V(\Gplus n)\cap V(\Ghor n)$. 
Let $w_iv$ be the last of all such edges on the path $P$,
starting from $w_F$. That is, let $v$ be such a vertex
that the path $P$ has the form $P=w_F\ldots w_iv \ldots u_1$ and its subpath $Q=v\ldots u_1$ 
contains no vertices from the set $\{w_0,w_1,\ldots,w_{F}\}$;
that is, $Q$ is a path in the original graph $\Graph n.$
We have $v\in V(\Ghor n).$ On the other hand, by assumption, $u_1\in V(\Gvert n)$.
%
Therefore, if $m$ is the length of $Q$ then by
 Theorem~\ref{v:Veta1}
 we have 
$m > l - 2 = F-1.$
So $(F- d(u_1))\varepsilon\leq (F- m) \varepsilon \leq 0$ and by the definition~(\ref{e:ThirdPartDefG}) of $\g n$
we have
$\g n(u_1)=0.$

Point {\bf{2}} of the statement follows from (\ref{e:WhereIsG}).
\end{proof}

\begin{theo}\label{v:Veta3}
There exist 
  functions $g,h\in C(\mathbb R)$ such that
  
\begin{enumerate}  
\item $|f(x)-g\circ p(x)-h\circ q(x)|\leq 6\varepsilon$
for all $x\in K$ and
\item $||g||\leq ||f||$  and $||h||\leq 2||f||.$
\end{enumerate}
  
\end{theo}

\begin{proof}
Let us choose $n$ in the same way as in the proof of Theorem~\ref{v:Veta2}.
First we define  functions $g$  and $h$ on $p(V(\Graph n))$ and  on
$q(V(\Graph n))$, respectively, and then we extend them to $\mathbb
R$.

For each point $x\in p(V(\Graph n))$ fix an arbitrary vertex $u\in
V(\Graph n)$ which lies in the $p$-fiber of $x$, i.e. such that $p(u)=x$, and define
 $$g(x)=\g n(u).$$ Similarly for each point $y\in q(V(\Graph
n))$ fix an arbitrary vertex $u\in V(\Graph n)$ such that $q(u)=y$
and define $$h(y)=\f n(u)- \g n(u).$$
Let $u\in  V(\Graph n)$ be an arbitrary vertex. Let us show that 
 $|\f n(u)- g\circ p(u) - h\circ q (u)|\leq 3\varepsilon$.
 %
Let $u_1,u_2 \in  V(\Graph n)$ be the vertices  such that
$p(u_1)=p(u)$, $q(u_2)=q(u)$ and  $g\circ p(u) = \g n(u_1)$,  $h \circ
q(u) = \f n(u_2) - \g n(u_2).$
%
%

%
If $|u-u_1|<\delta$, then  by condition
{\bf 1(a)} of Theorem~\ref{v:Veta2} we have $|\g n (u) - \g
n(u_1)|\leq\varepsilon$ and by Lemma~\ref{v:ZmenaFNaKratkej} we have $|\f n (u) - \f
n(u_1)|\leq\varepsilon$.  If  $|u-u_1|\geq \delta$, then
 by condition {\bf 1(c)} of the Theorem~\ref{v:Veta2}
 we have $\g n (u) = \g n(u_1)=0$.

Similarly, if $|u-u_2|<\delta$ then
we have 
 $|\g n (u) - \g
n(u_2)|\leq \varepsilon$ and  $|\f n (u) - \f
n(u_2)|\leq \varepsilon$, otherwise  by condition {\bf 1(b)} we have both
 $|\g n (u)-\f n(u)|\leq \varepsilon$ and
 $|\g n (u_2)-\f n(u_2)|\leq \varepsilon$.

By considering all four combinations of these possibilities we 
see that indeed
\begin{equation}\label{e:Veta2_DifferenceOnVertices}
 |\f n(u)- g\circ p(u) - h\circ q (u)| \leq 3\varepsilon.\end{equation}
%
%
%
%

Let $u,v\in V(\Graph n)$ be two vertices, which are
neighbors in the $p$-projection of the lattice,
i.e. such that $p(u)+
{1}/{2^n} = p(v).$ 
Again, there are vertices $u_1,v_1 \in  V(\Graph n)$  such that
$p(u_1)=p(u)$, $p(v_1)=p(v)$  and  $g\circ p(u) = \g n(u_1)$,  
$g\circ p(v) = \g n(v_1)$. Since $u_1$ and $v_1$ are neighbors in the vertical
projection of the lattice as well, $u_1v_1$ is a vertical edge in  $E(\Graph n)$.
 If is is a short edge, then by condition {\bf 1(a)}
of Theorem~\ref{v:Veta2} we have $|\g n(u_1) - \g n(v_1)|\leq \varepsilon$
and if it is a long edge, then by condition {\bf 1(c)} we have
$|\g n(u_1) - \g n(v_1)|=0.$ So in both cases
%
%
 \begin{equation}\label{e:Veta2_Difference OfG}|g\circ p (u) - g\circ p (v)|\leq
\varepsilon.\end{equation}
Similarly, if $u,v\in V(\Graph n)$ are such that $q(u)+ 1/2^n
= q(v) $ then
\begin{equation}\label{e:Veta2_Difference OfH}|h\circ q (u) - h\circ q
(v)|\leq 2\varepsilon .\end{equation}
%
%
%
%
We extend  $g$ and $h$ to $\mathbb R$ in the following way.   Let
us denote the vertical projection $p(V(\Graph n))$ of the vertex set by  
$\{x_1,x_2,\ldots,x_m\}$ with $x_1<x_2<\ldots<x_m$; by the definition
of $\Graph n$ we have $|x_i-x_{i+1}|\geq 1/2^n$ for all $n$.

For each $i$ such that $ x_i + 1/2^n=x_{i+1}$, we extend $g$ linearly
on $[x_i;x_{i+1}]$ (between the values $g(x_i)$ and $g(x_{i+1})$).

Now, let $i$ be such that $x_i + 1/2^n< x_{i+1}$.
On $[x_i;x_i + 1/2^n]$ we extend $g$ as a constant, equal to $g(x_i)$, and
on $[x_i + 1/2^n; x_{i+1}]$ we extend $g$ linearly. 
%
Note, that the definition of $\Graph n$ implies that 
 $p(K)\cap [x_i + 1/2^n; x_{i+1}) = \emptyset$, so we do not worry about
 the values of $g$ on interval $[x_i + 1/2^n; x_{i+1})$.

On $(-\infty;x_1]$ and $[x_m;\infty)$ we extend $g$ as a constant.


We extend the function $h$  to $\mathbb R$ similarly.

Hence $g$ and $h$ are continuous functions defined on $\mathbb R$
and $||g||\leq ||\g n||$, $||h||\leq ||\f n|| +
||\g n||.$ By condition {\bf 2} of Theorem~\ref{v:Veta2}
we have $||\g n|| \leq ||\f n||$ and from the construction of $\f n$ we have
$||\f n||\leq ||f||$. It follows that {\bf 2} of our statement holds true.


%
Let $x$ be an arbitrary point from the set $K$. Let $u \in V(\Graph n)$ be a vertex such that $x$ lies in the square
$\left[p(u);p(u) + {1}/{2^n}\right) \times \left[q(u);q(u) +
{1}/{2^n}\right).$
%
%
%
  The
definition of $g$ and (\ref{e:Veta2_Difference OfG}) imply that
$|g\circ p(u)-g\circ p(x)|\leq \varepsilon.$ Similarly the definition of $h$ and
(\ref{e:Veta2_Difference OfH}) imply that $|h\circ q(u)-h\circ
q(x)|\leq 2\varepsilon.$ The choice of $n$ and~(\ref{e:Veta2ChoiceDelta}) imply that $|\f n(u) -
f(x)|\leq \varepsilon$. Using  
(\ref{e:Veta2_DifferenceOnVertices}) we finally obtain
$$ |f(x)-g\circ p(x)-h\circ q(x)| \leq  6\varepsilon.$$

\end{proof}

The following result is a direct consequence of Theorem 4.13, implication (b) 
$\Rightarrow$ (c) of~\cite{rudin}.

\begin{theo}\label{v:Veta4} Let $K\subseteq \mathbb R^2$ be a compact subset of the plane.
If there exists $k\in \mathbb N$ such that for each function $f\in C(K)$
and each $\varepsilon>0$ there exist functions $g',h'\in C(\mathbb R)$
with
\begin{enumerate}
\item $|f(x)-g'\circ p(x)-h'\circ q(x)|\leq \varepsilon$
for all $x\in K$ and
\item $||g'||\leq k||f|| \text{\ \ and\ \ } ||h'||\leq k||f||$
\end{enumerate}
then there exist functions $g,h\in C(\mathbb R)$ such that
$$f(x)=g\circ p(x)+ h\circ q(x)$$
for all $x\in K$.
\end{theo}

\begin{proof}[Proof of the main Theorem~\ref{v:Basic}]
Theorem~\ref{v:Veta3} shows, that if a set $K$ contains no array on three points then the assumptions of
Theorem~\ref{v:Veta4}  are satisfied. Hence, the main result follows.
\end{proof}

\section{Alternative Proof}

In this section we give the proof of Theorem~\ref{v:Basic} suggested by an anonymous referee.
Let a compactum $K\subseteq \mathbb R^2$ contain no array on three points.
Then the fibers of $p$ and $q$ form an upper semi-countinuous decomposition of the set $K$.
Therefore we may consider the projection $\fun \psi K L$ to the factor space $L$.
Let $f\in C(K).$
Let the {\em variation\/} of $f$ on a subset $A$ of $K$ be the diameter of $f(A)$.
Fix $\varepsilon>0$ and denote by $F_p$ and $F_q$ the union of the fibers of $p$ and $q$,
respectively, on which the variation of $f$ is greater or equal than $\varepsilon$.
Let $\fun {s_p,s_q} L {[0;1]}$ be a partition of unity subordinated to the cover of $L$ consisting of the sets
$L \setminus \psi(F_p)$ and $L \setminus \psi(F_q)$, respectively.
Define functions ${f_p,f_q} \in C(K)$ by letting
$f_p(x)=f(x)\cdot s_p(\psi(x))$ and $f_q(x)=f(x)\cdot s_q( \psi(x))$.
Then the variations of $f_p$ on the fibers of $p$ and of $f_q$ on the fibers of $q$ are 
less than $\varepsilon$ and therefore there exist functions $g,h\in C(\mathbb R)$
such that $|f_p(x) - g\circ p(x)|<\varepsilon$ and $|f_q(x) - h \circ q(x)|<\varepsilon$ for all $x\in K$.
Thus $|f(x)-g\circ p(x) - h\circ q(x)|<2\varepsilon$ for all $x\in K$. The result follows from
Theorem~\ref{v:Veta4}.


\bibliographystyle{alpha}
\bibliography{Bibl_Trenklerova}

\end{document}